\documentclass[11pt]{article}

\usepackage{amsfonts}
\usepackage{amsmath,amsthm,amscd,amssymb,mathrsfs,setspace, textcomp,bm}
\usepackage{enumerate}
\usepackage{latexsym,epsf,epsfig}
\usepackage{color}
\usepackage[hmargin=2.25cm,vmargin=2.75cm]{geometry}
\usepackage{wrapfig}

\usepackage{cite}
\usepackage{hyperref}

\newcommand{\ds}{\displaystyle}

\newcommand{\xb}{{\bf{x}}}

\newcommand{\bu}{\mathbf u}

\newcommand{\h}{\hat}

\theoremstyle{plain}

\newtheorem{assumption}{Assumption}
\theoremstyle{remark}
\newtheorem{remark}{Remark}[section]
\numberwithin{equation}{section}
\numberwithin{theorem}{section}
\numberwithin{remark}{section}
\numberwithin{assumption}{section}
\numberwithin{condition}{section}

\begin{document}
\title{Dynamic Equations of Motion for Inextensible Beams and Plates}

 \author{\normalsize \begin{tabular}[t]{c@{\extracolsep{.6em}}c@{\extracolsep{.6em}}c@{\extracolsep{.6em}}c}
      Maria Deliyianni\footnote{University of Maryland, Baltimore County, MD}&   Kevin McHugh \footnote{Air Force Research Lab, Wright-Patterson AFB, OH}&  Justin T. Webster\footnote{University of Maryland, Baltimore County, MD} & Earl Dowell \footnote{Duke University}\\
\it mdeliy1@umbc.edu  ~~ &  \it  ~~kevin.mchugh.3@us.af.mil  & ~~\it websterj@umbc.edu & \it  ~~earl.dowell@duke.edu ~~
\end{tabular}}

\maketitle

\begin{abstract} {\noindent The large deflections of cantilevered beams and plates are modeled and discussed. Traditional nonlinear elastic models (e.g., that of von Karman) employ elastic restoring forces based on the effect of stretching on bending, and these are less applicable to cantilevers. Recent experimental work indicates that elastic cantilevers are subject to nonlinear inertial and stiffness effects. We review a recently established (quasilinear and nonlocal) cantilevered beam model, and consider some natural extensions to two dimensions---namely, inextensible plates. Our principal configuration is that of a thin, isotropic, homogeneous rectangular plate, clamped on one edge and free on the remaining three. We proceed through the geometric and elastic modeling to obtain equations of motion via Hamilton's principle for the appropriately specified energies. We enforce {\em effective} inextensibility constraints through Lagrange multipliers. Multiple plate analogs of the established 1D model are obtained, based on various assumptions. For each plate model, we present the modeling hypotheses and the resulting equations of motion. It total, we present three distinct nonlinear partial differential equation models, and, additionally, describe a class of ``higher order" models. Each model has particular advantages and drawbacks for both mathematical and engineering analyses. We conclude with an in depth discussion and comparison of the various systems and some analytical problems.
  \\[.15cm]
\noindent {\bf Key terms}:  cantilever, nonlinear elasticity, inextensibility, quasilinear PDE, energy harvesting
 \\[.15cm]
}
\\[.15cm] {\bf Distribution Statement}:  A. Approved for public release: distribution unlimited. AFRL-2021-2256.
\end{abstract}
\maketitle

\section{Introduction}

The purpose of this manuscript is to obtain equations of motion---including boundary conditions---for the large deflections of cantilevers. Cantilevers, as a class of elastic objects, have received less attention in the voluminous literature on beams and plates. In particular, there is a dearth of work analyzing the {\em large deflections} of cantilevers. Recent applications, as described below, have brought about the need for viable partial differential equation (PDE) models (i.e., distributed parameter systems) that capture the dynamics. For instance, we are interested in the dynamics of a long, slender flapping beam. In terms of spatial distribution,  we consider both beams (1D) and  plates (2D) here; as we are interested in large deflections,  models will be {\em nonlinear} in the displacement variables. The linear theory of elasticity is well established in both engineering and mathematics for beams and plates (e.g., \cite{beams,lagleug,novo,bolotin,ciarlet,antman}). Often, standard models---such as that of Euler-Bernoulli, Rayleigh, or Timoshenko for beams, and Kirchhoff-Love or von Karman for plates---are utilized because they are sufficiently accurate, while being analytically and computationally tractable. For cantilevered beams, we aim to produce a theory which generalizes the traditional linear cantilevered theory in the realm of (nonlinear) large deflections. We will consider {\em planar} transverse and in-axis displacements for beams, and for plates, we will consider rectangular plates, with 3D displacements from an equilibrium given in a standard Euclidean frame.

Unlike traditional nonlinear theories of beams and plates, where the boundaries are fully restricted, we will not assume that nonlinear forces result from local stretching (extensibility), as would be the case for fully clamped or hinged structures (panels). Rather, we will assume certain {\em inextensibility} conditions---that the beam does not stretch in specified ways throughout deflection. The enforcement thereof via Lagrange multipliers in the dynamics will yield nonlinear inertial terms in the equations of motion. Additionally, for consistency in {\em order considerations}, we will invoke high order expressions for the potential energy, themselves simplified via inextensibility constraints. The additional contribution will be nonlinear stiffness terms in the equations of motion when Hamilton's principle is invoked.

\subsection{Motivating Application}
Under harmonic excitation \cite{dowell3}, follower forces \cite{follower}, or other boundary forces, cantilevers may experience large deflections (say, on the order of their length). We focus here on the motivating application of {\em piezoelectric energy harvesters}. Slender cantilevers in axial\footnote{along the principal length of the beam or plate} airflows can experience a self-destabilization known as {\em flutter}, even at markedly low flow velocities \cite{dowell4,tang,ELSS,gibbs}. Beyond some critical flow velocity, the flow-structure system enters a {\em limit cycle oscillation} (LCO), and such dynamics can generate extractable power \cite{DOWELL}.
\begin{wrapfigure}{R}{3.25in} 
\vspace{-20pt}
  \begin{center}
    \includegraphics[width=3in]{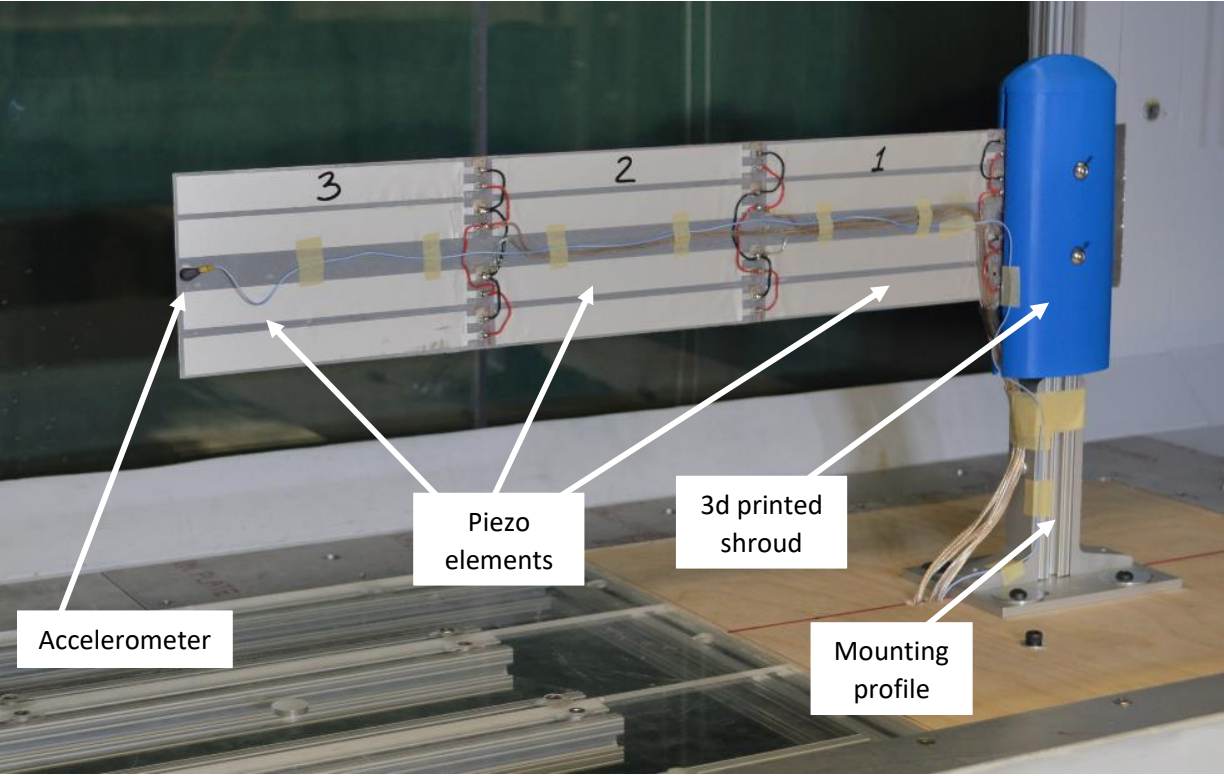}
    \caption{\small Experimental setup from \cite{experimental}.}
    \label{fig:piezo}
  \end{center}
  \vspace{-17pt}
  \vspace{2pt}
\end{wrapfigure} 

While interest in airfoil and panel flutter has been {\em immense} for 75 years (see \cite{dowell1} and references therein), interest in the motion of cantilevers driven by axial flow has been modest, at least until about 20 years ago. In \cite{sleepapnea2,huang}, the cantilever flutter of the human soft palate was considered, and, as in most engineering scenarios, the flutter thereof is undesirable, with a design goal of prevention. Lack of  interest for this configuration seems to derive from the few practical motivations. However, recent interest in alternative energies, in particular, energy harvesting \cite{energyharvesting,DOWELL}, have provided reasons for the investigation of cantilever LCOs. In this application, dynamic instability is {\em encouraged} in order to extract energy from post-onset cantilever LCOs. This has been shown possible for {flow-induced cantilever oscillations} with affixed piezoelectric material \cite{piezomass,energyharvesting,shubov,TP1}.

This motivates the problem of characterizing  {\em post-flutter} behaviors of cantilever (i.e., LCOs), which necessitates a viable set of nonlinear equations of motion. In reality, such a structural analysis should be truly 2D, as we must practically consider distributed and variable stiffness and damping effects, owing to the use and placement of piezoelectric patches and layers. Recent experiments \cite{dowellmaterial,experimental,dowelleffects} and others \cite{energyharvesting,piezomass} have studied piezoelectric structures in elastic energy harvesters. 

\vskip.2cm
\noindent{\em Thus, there is a need for a {theory} of large deflections for cantilevers. In physically applicable scenarios, the relevant theory should be truly 2D. To date, very little work has been done---experimentally and analytically---in this configuration, relative to the immense body of work on {extensible} nonlinear plates}.

\subsection{Modeling and Analysis of Cantilever Large Deflections}

The primary focus of this work is to produce nonlinear equations of motion called for above for 1D and 2D cantilevers undergoing large deflections. These equations can subsequently be analyzed numerically and analytically. As we have established, the models of interest are based on inextensibility, rather than extensibility.  In the case of extensible structures, models and theory are well-established, especially in scenarios where the boundary is fully restricted (e.g., paneling \cite{bolotin,springer,dowell1}); such nonlinear restoring forces are based on the effect of local stretching on bending. 

In the case of beams, a prominent example is the Woinowsky-Krieger beam \cite{wonkrieg} which appends a semilinear extensible force to any of the standard linear beams such as those of Euler-Bernoulli, Rayleigh, shear, or Timoshenko \cite{beams}. We note some classical mathematical references for the analysis of such extensible beams, e.g., \cite{ball,edenmil}, as well as some more recent \cite{htw,hhww} that specifically address extensible cantilever systems. Alternatively, \cite{lagleug} provides a modeling account (with subsequent analysis of solutions and their stability) for an extensible beam that is intimately related to the so called {\em full von Karman} model for plates. 

Indeed, for plates, the linear theory of Kirchhoff-Love is well-established, and the prominent nonlinear generalization for large deflections is the theory of von Karman (\cite{ciarlet} and references therein). One can consult \cite{bolotin,novo,berger0} for classical discussions, and the rigorous mathematical justifications of these models in \cite{ciarlet,laglio}. The more modern references of \cite{lagnese,springer} provide detailed discussion of solutions, stability, and dynamical systems aspects of von Karman dynamics (among others). We finally mention \cite{koch} for the analysis of the {\em full von Karman} plate, which is an extensible plate system that accounts for both in-plane and transverse inertial terms (in contrast to the scalar von Karman equations, such as those in thoroughly analyzed in \cite{springer}); the full von Karman plate equations, as we shall see, have some mathematical similarities to the plate models derived herein.

We forgo an extensive discussion of the literature on the analysis of the aforementioned nonlinear beams and plate models, and focus specifically here on recent work that {\em attempts to capture the effects of inextensibility in cantilever dynamics}. It should be noted that the theory of von Karman is not viable (namely, the operative hypotheses which invoke the effect of stretching on bending) in the case of a cantilever. As such, a theory based on inextensibility must be developed from first principles. The earlier work of Pa\"{i}doussis derives and discusses 1D inextensible pipes conveying fluid, along with associated numerical stability analyses \cite{semler2, paidoussis}. These works derive the PDE equations of motion and boundary conditions for 1D structures. Subsequently, the thesis \cite{techthesis} focuses on inextensible beams and plates, but primarily focuses on finite dimensional descriptions through Rayleigh-Ritz considerations, not providing the explicit equations of motion in the standard Euclidean frame. The work of Dowell et al. further elaborated on inextensible beams and plates, developing several related approaches and investigating an inextensible theory experimentally and numerically in the sequence of papers \cite{inext2, tang, dowell4}. Finally, the work of McHugh et al. \cite{inext1,kevinthesis,follower} provided the first explicit calculus of variations derivation of the PDE equations of motion of an inextensible beam (in both the free-free and cantilever configurations).  Apart from these papers, other aforementioned analyses primarily employ the Rayleigh-Ritz procedure to obtain finite dimensional equations of motion directly from specified energies, circumventing the PDE equations of motion. Numerical comparisons and computational studies are then performed. The discussion of boundary conditions for the 2D system (from first principles) is typically omitted, and the rigorous analysis of solutions, as well as discussions of well-posedness thereof, are not made. 

Recently, for the intextensible cantilever as derived in \cite{inext1}, a mathematical theory of solutions has been developed \cite{DW,DGHW}. In the present manuscript, we aim to take an analogous first step in this direction, by providing the PDE equations of motion for the ``natural" beam extension to a 2D inextensible plate. We specifically follow the approaches for the beam, coming from references \cite{inext1,kevinthesis} and \cite{DGHW}. To begin our considerations, we will first revisit the setup for an inextensible cantilevered beam; this will provide a template for plate considerations. In doing so, we will include a novel ``higher order" beam model, which is motivated from a possible need for inclusion of such effects for the plate.
 
 As we shall see, there is no single, clear plate extension of the inextensible, cantilevered beam model. Indeed, there are certain modeling and order choices, each of which yields a different system for the equations of motion. The presence of 2D shear effects is a critical aspect of the analysis of plates here. Providing a careful derivation of the equations of motion is indeed quite necessary, as we will see, since  ``natural" modeling choices produce nontrival nonlinear boundary conditions for each plate model. {\em This is to say: it is notable that for the 2D plate, the modeling hypotheses produce (nonlinear) boundary conditions which are {not} those of a standard, linear, cantilevered plate. Moreover, some aesthetic simplifications which present themselves in 1D will be conspicuously absent in 2D.}

 In the end, we will present various sets of hypotheses, yielding a variety of different systems of equations of motion. {\em To the best knowledge of the authors, each of the plate systems here is novel}. In future work, we will numerically analyze these models to compare against \cite{dowell4,tang3}. Moreover, the theory developed in \cite{DW} will be adapted, at least in the cases where the 2D equations of motion most closely resemble those of the inextensible cantilevered beam in \cite{inext1,DW}. Indeed, one principal goal in this work is to address inextensibility in the standard coordinate framework---which is to say, to extend standard linear beam theory, rather than using more sophisticated geometric models (such as those of \cite{antman}). In this way, existing mathematical and numerical tools, as well as theory, can be utilized, viewing these models as ``extensions" of the classical linear theory of elasticity. Finally, we provide a brief discussion of the numerical work from \cite{kevinthesis} in Section \ref{Discussion}.
	
\subsection{Outline of the Remainder of the Paper}
To provide context for the approach we take to {\em inextensible plates}, we first give a summary of the work in \cite{inext1,DW} for the recent {\em inextensible cantilevered beams} in Section \ref{Beams}. 

Section \ref{Plates} is the principal modeling section. There, we walk through our choices for the potential energies, as well as the interpretation and enforcement of ``inextensibility" in the plate. Ultimately, we present three distinct models in Section \ref{Plates}, with clear discussion of the relevant modeling and order hypotheses. Each model is presented as a system of partial differential equations  in the elastic deflection variables, as well as constraint variables. We conclude Section \ref{Plates} with a description of how to obtain a class of higher order models, though without explicitly producing the equations of motion.

Finally, in Section \ref{Discussion} we discuss and compare the derived models from both the engineering and analysis-of-PDEs points of view. 

\subsection{Conventions and Notation}
We utilize the convention that spatial points $\mathbf x$ are identified with their position vector $\langle x,y,z\rangle$. We will sometimes write $x_1=x, ~x_2=y,~x_3=z$ when it is notationally expedient. The equilibrium (undeformed) domain under consideration will be cylindrical, namely, $\{\mathbf x \in \Omega \times (-h/2,h/2)\}.$ The domain  $\Omega$ will be identified as the centerline, $z=x_3=0$. We make use of the partial derivative notation $\dfrac{\partial}{\partial x_i}=\partial_{x_i}$, as well as the standard 2D spatial gradient $\nabla = \langle \partial_{x},\partial_{y}\rangle$, divergence $[\nabla \cdot]$, and Laplacian (in rectangular coordinates) $\ds \Delta = \nabla \cdot \nabla = \sum_{i=1}^2 \partial_{x_i}^2$.


\section{Inextensible Cantilevered Beams}\label{Beams}
We begin our discussion with cantilevered beams. These are markedly simpler, owing to the lack of chord-wise and span-wise interaction. The exposition here mirrors \cite{inext1,DGHW}. By presenting a streamlined version of the 1D model, and its procurement from precise geometric and order hypotheses, we elucidate the decisions which will lead to the variety of plate models in the subsequent sections. 

Consider a slender, isotropic, homogeneous beam, as is standard in Euler-Bernoulli beam theory \cite{beams}.  We assume that the filaments of the beam remain perpendicular to the centerline throughout deflection (the Kirchhoff-Love hypothesis---see \cite{beams,lagleug}). 
The variable $x \in [0,L]$ will represent the beam's centerline at equilibrium. We utilized the notation $\langle x+u(x,t) , w(x,t)\rangle \in \mathbb R^2$ for the planar deflection that corresponds to the equilibrium point $x$ at instant $t$. This is to say that $u$ is the in-axis displacement of the beam, whereas $w$ is the transverse displacement. In this way, the beam's displaced curve is parametrized by $x$ with position vector $\langle x+u(x,t) , w(x,t)\rangle$ at time $t$.
If we let $\varepsilon$ denote the strain associated to the beam centerline, we have the identity \cite{lagleug,novo,Culver}:
\begin{equation}\label{strainexpression1}\big[1+\varepsilon]^2=(1+u_x)^2+w_x^2.\end{equation}
\subsection{Inextensibility}\label{inextbeamsec}
When the beam is {\em inextensible}, we take there to be no center-line extensional stress, i.e., $\varepsilon(x,t)=0$. This  yields the condition:
\begin{equation}\label{realinext}
1=(1+u_x)^2+w_x^2.
\end{equation}
By expanding \eqref{realinext}, we see  that if $w_x \sim \eta$, we will have $u_x \sim \eta^2$:
$$2u_x+u_x^2+w_x^2=0.$$
Therefore, we may elect to drop terms of order $u_x^2 \sim \eta^4$, owing to their relative order. Approximating, then
 $$0=2u_x+w_x^2 ~~\Rightarrow~~ u_x=-\frac{1}{2}w_x^2.$$ The above can be taken an {\em effective inextensibility constraint}, providing a direct relationship between $u$ and $w$: \begin{equation}\label{inext} u(x,t)=-\frac{1}{2}\int_0^x[w_x(\xi,t)]^2d\xi.\end{equation}

 We note that the above analysis can easily be extended to ``higher order", as follows: 
\begin{equation}\label{realinext*}
1=(1+u_x)^2+w_x^2~~ \implies~~ u_x=-1\pm\sqrt{1-w_x^2}.
\end{equation}
The binomial expansion can then be used to obtain expressions to higher order
$$ u_x=-\frac{1}{2}w_x^2-\frac{1}{8}w_x^4-\frac{1}{16}w_x^6...,~~|w_x|<1.$$

Let us now fix the nomenclature  for {effective inextensibility constraints}, with reference to the order of approximation.
\begin{assumption} In referring to $\eta^2$-order and $\eta^4$-order approximations for the inextensible beam, we respectively assume one of the following relations:
\begin{equation}\begin{cases}
\text{$\eta^2$ order}:&\ds ~~u_x=-\frac{1}{2}w_x^2   \\[.4cm] 
 \text{$\eta^4$ order:~~}&\ds ~~u_x=-\frac{1}{2}w_x^2-\frac{1}{8}w_x^4. 
\end{cases}
\end{equation}
\end{assumption}

\subsection{Energies}

Define the elastic potential energy ($E_P$) via beam curvature $\kappa$ and constant stiffness $D$ \cite{bolotin,semler2,huang,novo}\footnote{beam flexural rigidity, which can be given in terms of inertial and Young's coefficients, $EI$} in the standard way:
 $$E_P \equiv \frac{D}{2}\int_0^L\kappa^2dx.$$ 
As discussed above, the beam's displaced configuration is parametrized by $x \in [0,L]$. As such, the curvature is given by \cite{Culver} $$\kappa = \dfrac{(1+u_x)w_{xx}+u_{xx}w_x}{[(1+u_x)^2+w_x^2]^{3/2}}.$$ 
 Substituting in the (full) inextensibility constraint \eqref{realinext} in to the curvature $\kappa$, we obtain:
\begin{align*}
\kappa = (1+u_{x})w_{xx} - w_{x}u_{xx} =&~ (1-w^2_{x})^{1/2}w_{xx} - w_{x}(w_{x}w_{xx}(1-w^2_{x})^{-1/2}) =  \frac{w_{xx}}{(1-w^2_{x})^{1/2}} .
\end{align*}
In $E_P$, we may invoke the appropriate geometric series for $\kappa^2=\dfrac{w_{xx}^2}{(1-w_x^2)}$ and truncate to a given order {\em in the nonlinear coefficient} of $w_{xx}$; we obtain   ~$\ds \kappa^2 \approx w_{xx}^2(1+w_x^2)$ in the case of the $\eta^2$-order approximation, and $\ds \kappa^2 \approx w_{xx}^2(1+w_x^2+w_x^4)$ in the case of the $\eta^4$-order approximation. \begin{remark} Note that, at the stage of the potential energy simplification, we attempt to {\em match} the order of approximation for the inextensibility condition and the approximated potential energy. \end{remark}

 The kinetic energy ($E_K$) is defined in the standard way, where we have mass-normalized the beam
  $$E_K \equiv \frac{1}{2}\int_0^L \left[ w^2_t + u^2_t \right]dx.$$ 
\subsection{Beam Equations of Motion for $\eta^2$-order}\label{beamsec1}
We consider displacements $u$ and $w$ (and virtual displacements $\delta u$ and $\delta w$) which respect the essential boundary conditions at $x=0$:
\begin{align*} &w,~w_x,~\delta w, ~\delta w_x:~ 0 ~\text{ at }~x=0;&~~~u,~\delta u:~0 ~\text{ at }~x=0.& \end{align*}
The effective inextensibility constraint, defined as $f \equiv u_{x} + (1/2)w^2_{x}=0$, is enforced through the {\em Lagrange multiplier $\lambda$}. We will minimize the Lagrangian:
\begin{equation}
\label{Lagrangian}
\mathcal{L} = E_{K} - E_{P} + \int_0^L \lambda f dx.
\end{equation}
Invoking Hamilton's principle for the time-integrated variation of $(\ref{Lagrangian})$ yields the Euler-Lagrange equations of motion and the associated boundary conditions (after the requisite integration by parts) for the displacements $u$ and $w$.
\begin{align}
\label{deltau}
\text{from arbitrariness of}~\delta u:& ~~~~u_{tt} + \lambda_x = 0;\\\label{deltaw}
~\text{from arbitrariness of}~\delta w:& ~~~~w_{tt} - D \partial_{x} \left( w^2_{xx} w_{x} \right) + D \partial_{xx} \left( w_{xx} \left [ 1 + w^2_{x} \right ] \right) + \partial_x \left( \lambda w_{x} \right) = 0.
\end{align}
For the (natural) boundary conditions at $x=L$ in $w$ (with $u$ and $\lambda$ then inferred) we have:
\begin{align}\label{firstnatural}
& \lambda(L) = 0;~~~
(1+w_x^2(L))w_{xx}(L) = 0 ; ~~~
(1+w_x^2(L))w_{xxx}(L)+w_x(L) w^2_{xx}(L)=0.&
\end{align}
The standard free beam boundary conditions---$w_{xx}(L)=w_{xxx}(L)=0$---
are then obtained algebraically from \eqref{firstnatural}, and from equation \eqref{deltau} we write $\lambda(L)-\lambda(x) = - \int_x^L u_{tt}(\xi) d\xi.$
We invoke the previously obtained $\lambda(L)=0$ to obtain ~$\ds
\lambda(x) =  \int_x^L u_{tt}(\xi) d\xi.$

Presenting the system in $(u,w,\lambda)$ yields some degree of redundancy, as one can fully eliminate both $\lambda$ and $u$ in describing the dynamics above. Indeed, we can solve for $\lambda$ in terms of $u$ from \eqref{deltau}, then we can recover $u$ in terms of $w$ from the  constraint \eqref{inext} to present the system solely in $w$. However, we present the system in its full form in anticipation of more complicated 2D systems to follow; in the models presented below, we will not be able to eliminate all non-transverse variables from the ``complete" description.
We can summarize, then, in the following description of the {\em $\eta^2$ inextensible beam} system:
\begin{align}
u_{tt} + &~ \lambda_x = 0\\
w_{tt} -& ~D \partial_{x} \left( w^2_{xx} w_{x} \right) + D \partial_{xx} \left( w_{xx} \left [ 1 + w^2_{x} \right ] \right) + \partial_x \left( \lambda w_{x} \right) = 0\\ \label{inextuse}
u_x=&-\frac{1}{2}w_x^2\\
w(0)&=w_x(0)=~0;~w_{xx}(L)=w_{xxx}(L)=0; \\
u(0)&=0;~~u_x(L)=-\frac{1}{2}w_x^2(L) \\
\lambda(0)& = \int_0^L u_{tt}(\xi) d \xi;~~\lambda(L)=0.
\end{align}
Of course, the above equations would be supplemented by appropriate initial conditions for $w$, namely:
$$w(x,0)=w_0(x),~~w_t(x,0)=w_1(x).$$
Initial conditions for $u$ are determined from those for $w$ through the inextensibility constraint.

\subsection{Beam Equations of Motion to $\eta^4$-order}\label{quarticbeam}
Repeating the steps in the previous the section, we retain terms up to $\eta^4$ in the effective inextensibility condition and the curvature expression $\kappa^2$. We obtain---paying particular attention to the natural boundary conditions---the {\em $\eta^4$ inextensible beam} system (omitting initial conditions):
\begin{align}
u_{tt} +&~ \lambda_x = 0\\
w_{tt} -& ~D \partial_{x} \left( [w_x+2w_x^3]w^2_{xx} \right) + D \partial_{x}^2 \left( w_{xx} \left [ 1 + w^2_{x}+w_x^4 \right ] \right)+ \partial_x \left( \lambda w_{x} \right) = 0\\
u_x=&~-\frac{1}{2}w_x^2-\frac{1}{8}w_x^4\\
w(0)&=w_x(0)=0;~w_{xx}(L)=w_{xxx}(L)=0;\\
u(0)&=0;~u_x(L)=-\frac{1}{2}w_x^2(L)-\frac{1}{8}w_x^4(L) \\
\lambda(0)& = \int_0^L u_{tt}(\xi) d \xi;~~\lambda(L)=0.
\end{align}

\section{Inextensible Cantilevered Plates}\label{Plates}

Let $\mathbf x = \langle x, y \rangle \in \mathbb R^2$. 
We will let $\langle u(\xb,t), v(\xb,t) ,w(\xb,t) \rangle \in \mathbb R^3$ denote the mid-plate displacement (from equilibrium) of a rectangular, cantilevered plate that occupies (at equilibrium) the region $\Omega \times (-h/2,h/2)$. The coordinate $u$ is span-wise displacement, with $v$ being chord-wise; $w$ is the transverse deflection. (Thus, $\langle x+u(\xb,t), y+v(\xb,t) ,w(\xb,t) \rangle$ describes the position of the equilibrium point $\mathbf x$ at the instant $t$.) The physical quantity $\nu \in (0,1/2)$ represents the so called {\em Poisson Ratio}.

Consider  the open rectangle $\Omega \equiv \{(x,y) \in (0,L_x) \times (0,L_y)\}$, representing the undeformed mid-plane ($z=0$) of the thin, homogeneous, isotropic plate. Let us take the four components of the boundary $\partial \Omega = \Gamma$ (in standard orientation) to be given by 
\begin{align*}
\Gamma_E = ~\{x=L_x,~y \in (0,L_y)\}, &&
 \Gamma_N=\{x \in (0,L_x), ~y=L_y\}, \\
 \Gamma_W=\{x=0,~y \in (0,L_y)\},&&
  \Gamma_S= ~\{x \in (0,L_x),~y =0\}.
  \end{align*}
\begin{center}
\includegraphics[scale=.8]{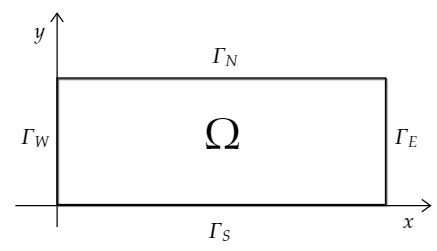}
\end{center}
We have the associated outward unit normals: ~$\mathbf n_E=\mathbf e_1,~ \mathbf n_N=\mathbf e_2,~ \mathbf n_W=-\mathbf e_1,~ \mathbf n_S=-\mathbf e_2$ and corresponding tangentials (respectively):~ $\mathbf e_2, -\mathbf e_1, -\mathbf e_2, \mathbf e_1$.

\begin{remark}\label{linBC} For reference, we provide the classical clamped-free conditions associated to a {\em linear} cantilever in this configuration, that is, the plate is clamped on $\Gamma_W$ and free elsewhere. 
On $\Gamma_N, \Gamma_S$ the boundary conditions are:
$$ \nu w_{xx}+w_{yy}=0;~~w_{yyy}+(2- \nu)w_{xxy}=0.$$
On $\Gamma_E$ the boundary conditions are:
$$w_{xx}+ \nu w_{yy}=0;~~w_{xxx}+(2- \nu)w_{yyx}=0.$$
And, of course, on $\Gamma_W$ we have the clamped conditions:
$$w=0;~~w_x=0.$$
The boundary conditions here are markedly simpler than the general case \cite{lagnese,springer} involving the well-known $B_1$ and $B_2$ operators, owing to the simplified rectangular geometry.
\end{remark}

\subsection{Full Elastic Potential Energy and Strain-Displacement Relations}\label{fullpotential}
In this section we discuss the plate potential energy. We have in mind to simplify it in accordance with the relevant plate inextensibility constraints, as done analogously for the beam.  We begin with the bulk description of the potential energy, namely that the potential energy $E_P$ has the form:
$$E_P=\frac{1}{2}\int_{-h/2}^{h/2}\int_{\Omega} \hat{\varepsilon}_{ij} \hat{\sigma}_{ij} d\Omega dz, $$
where $\hat \varepsilon$ and $\hat \sigma$ are the 3D strain and stress tensor (resp.) for 3D displacements. With the plate taken as homogeneous and isotropic,  we invoke the standard Hookean stress-strain relation
\begin{equation}\label{hooks}
\hat \sigma_{ij} = \frac{E}{1+\nu}\Big(\hat \varepsilon_{ij}+\frac{\nu}{1-2\nu}\hat \varepsilon_{kk}\delta_{ij}\Big),~i,j=1,2,3,
\end{equation}
where $E>0$ is the Young's modulus, $\nu \in (0,1/2)$ is the Poisson ratio, and the Einstein summation convention is in force (above $\delta_{ij}$ is that of Kroenecker).

As is customary in thin plate theory, we assume the Kirchhoff-Love hypotheses. 
\begin{assumption} Assume $\hat \sigma_{33}=0$. We also assume that plate filaments that are originally orthogonal to the mid-surface remain so throughout deflection, and that said filaments have constant length; this induces $\hat\varepsilon_{13}=\hat \varepsilon_{23}=0$. \end{assumption}
\noindent The above assumption reduces $E_P$ to the following expression \cite{lagnese,ciarlet}, which we will use throughout this treatment:
\begin{equation}
\label{potentialfoundamental}
E_P = \frac{1}{2} \frac{E}{1-\nu^2} \int_{\Omega}\int_{-h/2}^{h/2} \left[ \h \varepsilon^2_{11} + \h \varepsilon_{22}^2 + 2 \nu \h \varepsilon_{11} \h \varepsilon_{22} + \frac{1-\nu}{2}\h \varepsilon^2_{12} \right] dz d \Omega.
\end{equation}
We now attempt to simplify, utilizing only mid-plane strains $\varepsilon$, curvature $\kappa$, and the higher order tensor $\mu$ \cite{novo}, yielding:
\begin{align}\label{bulk1}
 \h \varepsilon_{11} & =  { \varepsilon} _{11} +z \kappa_{11} + z^2 \mu_{11}  \\[1em] \label{bulk2}
\h \varepsilon_{22} & =  { \varepsilon} _{22} +z \kappa_{22} + z^2 \mu_{22}  \\[1em] \label{bulk3}
 \h \varepsilon_{12} & =  { \varepsilon} _{12} +z \kappa_{12} + z^2 \mu_{12}. 
\end{align}
In all of what follows we drop the terms scaled by $z^2$ above, consistent with \cite{novo,kevinthesis}.
\begin{assumption}\label{z} In the expression for strains $\hat{\varepsilon}_{ij}$ for $i,j=1,2$, we assume that each $z^2\mu_{ij}$ for $i,j=1,2$, are negligible.\end{assumption}

It is necessary, at this point, to invoke strain-displacement relations \cite{novo}. Recalling that $u,v$ and $w$ denote the Lagrangian mid-plane displacements, we have
\begin{align}\label{strain1}
\varepsilon_{11} & = u_{x} + \frac{1}{2} \left[ u_x^2  + v^2_x + w_x^2 \right]\\[1em]
\varepsilon_{22} & = v_{y} + \frac{1}{2} \left[ u_y^2  + v^2_y + w_y^2 \right]\\[1em]\label{strain3}
\varepsilon_{12} & =u_y + v_x + u_x u_y + v_x v_y + w_x w_y,
\end{align}
which, as mentioned above, are the mid-plane extensible axial strains and the shear strain (resp.).

\subsection{Inextensibility}
The ``true" inextensibility conditions would correspond to taking ${\varepsilon}_{ij}=0$ \cite{novo,inext2}; this translates to no extensional stresses along the mid-plane of the plate. In that case, we have three conditions, taken from \eqref{strain1}--\eqref{strain3} by equating each to zero:
\begin{align}\label{fullinext1}
(1+u_x)^2+v_x^2+w_x^2=&~1 \\ \label{fullinext2}
u_y^2+(1+v_y)^2+w_y^2=&~1 \\ \label{fullinext3}
u_y+v_x+u_xu_y+v_xv_y+w_xw_y=&~0.
\end{align}
We refer to these as the {\em full} plate inextensibility conditions. The first condition is a span-wise inextensibility constraint, the second is chord-wise, with the lattermost corresponding to shear. 

The first two inextensibility conditions can be re-written (discarding non-physical square roots) as:
\begin{align}
u_x=&~-1+\sqrt{1-(v_x^2+w_x^2)} \\
v_y=&~-1+\sqrt{1-(u_y^2+w_y^2)}.
\end{align}
These can be truncated through Taylor expansions about equilibrium. Namely (up to quartic terms):
\begin{align}\label{taylor1}
u_x=-1+\sqrt{1-(v_x^2+w_x^2)} = & ~-\frac{1}{2}v_x^2-\frac{1}{2}w_x^2-\frac{1}{4}v_x^2w_x^2-\frac{1}{8}v_x^4-\frac{1}{8}w_x^4-... \\ \label{taylor2}
v_y=-1+\sqrt{1-(u_y^2+w_y^2)}=&~-\frac{1}{2}u_y^2-\frac{1}{2}w_y^2-\frac{1}{4}u_y^2w_y^2-\frac{1}{8}u_y^4-\frac{1}{8}w_y^4-...
\end{align}
\begin{remark} It is clear that any order analysis associated to truncating these identities---such as what is done for the beam in Section \ref{inextbeamsec}---will involve the quantities $v_x, u_y, w_x,$ and $w_y$. Specifically, the above invites comparisons of $~u_x$ to both $~v_x, w_x$ and $~v_y$ to both $~u_y, w_y$. \end{remark}

\subsubsection{Effective Inextensibility Conditions}
In working with the potential energy, we will invoke approximated versions of the inextensibility constraints \eqref{fullinext1}--\eqref{fullinext3}, analogously with the 1D beam. However, in the case of the plate, there are more complex order comparisons which can be made, and a variety of assumptions can be taken.

In the simplest case, we can assert the following $\eta^2$-type scaling assumption:

\begin{assumption}\label{scaling}Assume that $\partial_{x_i}w \sim \eta$ and assume that $[\partial_{x_i} u], [\partial_{x_i} v] \sim \eta^2$, for all $i=1,2$.\end{assumption}  This is to say that we assume the in-plane gradient is of higher order than the slopes associated to transverse deflection. With this assumption, we can reduce the full inextensibility constraints in \eqref{fullinext1}--\eqref{fullinext3} up to the order $\eta^2$ (thereby dropping terms of the form $[\partial_{x_i}u][\partial_{x_j}v]$, $[\partial_{x_i}u]^2$, $[\partial_{x_i}v]^2$). 
This produces a particular set of $\eta^2$-effective inextensibilty constraints:
\begin{align}\label{quadcon1}
u_x=& ~-\frac{1}{2}w_x^2 \\ \label{quadcon2}
v_y=& ~-\frac{1}{2}w_y^2 \\ 
u_y+v_x=&~-w_xw_y.\label{quadcon3}
\end{align}

\begin{remark} One can directly obtain, from the above equations, the identities:
$$4u_xv_y-2u_yv_x=[u_y^2+v_x^2]~~\text{ and }~~w_{xx}w_{yy}=w_{xy}^2.$$ 
These are elaborated upon in \cite{kevinthesis,inext2}. We note that neither of these ``composite" constraints over-constrain the model, since they are derived from (combinations of) prior equations that each eventually correspond (one-to-one) with a Lagrange multiplier variable.
 \end{remark}
 \begin{remark}
 The latter equality in the previous remark does provide another interesting identity via the boundary conditions and integration by parts:
 $$\int_{\Omega}w_{xy}^2d\mathbf x=\int_{\Omega}w_{xx}w_{yy}d\mathbf x = \int_{\Omega}w_{xy}^2d\mathbf x+w_{x}w_{yy}\Big|_{x=0}^{x=L_x}-w_xw_{xy}\Big|_{y=0}^{y=L_y}.$$ As we have yet to derive the higher order boundary conditions via a variational procedure, the only condition we can invoke is the essential condition that $w_x =0$ on $\{x=0\}$, which finally yields:
 $$w_{yy}(L_x)=w_xw_{xy}\Big|_{y=0}^{y=L_y}.$$
 \end{remark}

\begin{remark}
Alternatively, in the choice of truncation above in Assumption \ref{scaling}, we could more closely follow the logic of our beam approximation in simplifying \eqref{fullinext1}--\eqref{fullinext3}. Namely, we might assume only that $w_x\sim\eta$, from \eqref{fullinext1}, and $w_y \sim \eta$ from \eqref{fullinext2}. We could expand to see that $u_x\sim \eta^2$ with $v_x \sim \eta$, as well as $v_y \sim \eta^2$ with $u_y \sim \eta$. In this case, we can choose to retain terms in the Taylor expansions \eqref{taylor1}--\eqref{taylor2} up to order $\eta^4$:
\begin{align}
u_x=&  ~-\frac{1}{2}v_x^2-\frac{1}{2}w_x^2-\frac{1}{4}v_x^2w_x^2-\frac{1}{8}v_x^4-\frac{1}{8}w_x^4\\
v_y=&~-\frac{1}{2}u_y^2-\frac{1}{2}w_y^2-\frac{1}{4}u_y^2w_y^2-\frac{1}{8}u_y^4-\frac{1}{8}w_y^4\\
u_y+& v_x=-u_xu_y-v_xv_y-w_xw_y.
\end{align}
\end{remark}
Note that, if we go up to order $\eta^4$, the shear inextensibility condition will be retained in full---no truncation is called for. In the final class of models  in Section \ref{higherorder}, we will discuss the use of higher order approximations of the inextensibility constraints, analogous to the model in Section \ref{quarticbeam} 

\subsubsection{Curvature Expressions}\label{curvexp}
To obtain expressions for $\kappa_{ij}$ of the $(x,y) \in \Omega$ middle surface of the plate, we follow the Kirchhoff-Love hypotheses ($\hat{\varepsilon}_{13}=\hat{\varepsilon}_{23}=\hat{\varepsilon}_{33}=0$) and from \cite{novo}, we obtain\footnote{In \cite{novo}, the relations $\hat u = u+z\theta$, $\hat v = v +z\psi$, and $\hat w =w+ z\chi$ are plugged in to the strain relations for $\hat{\varepsilon}_{i3}$, then themselves simplified by geometric order considerations.}: 
\begin{align}
\kappa_{11} & =  (1+u_x)\theta_x+v_x\psi_x+w_x\chi_x  \\[.1cm]
\kappa_{22} & = u_y\theta_y+(1+v_y)\psi_y+w_y\chi_y  \\[.1cm]
\kappa_{12} & = (1+u_x) \theta_y+(1+v_y)\psi_x+u_y\theta_x+v_x\psi_y+w_x\chi_y+w_y\chi_x,
\end{align}
with $\theta, \psi, \chi$ given by:
\begin{align}
\label{curvaturenotsimplified*}
\theta  = -(1+v_y)w_x+v_xw_y;~~~~~
\psi  = -(1+u_x)w_y+u_yw_x;~~~~~
\chi  = u_x+v_y+u_xv_y-u_yv_x.
\end{align}
\begin{remark} As noted in \cite{novo}, the above expressions assume that strains are small compared to $1$. \end{remark}
We can thus obtain:
\begin{align} \label{curve1}
\kappa_{11} =&~(1+u_x)[v_{xx}w_y+v_xw_{xy}-v_{xy}w_x-(1+v_y)w_{xx} ]\\[.1cm]\nonumber
&+v_x[u_{xy}w_x+u_yw_{xx}-u_{xx}w_y-(1+u_x)w_{xy}]\\[.1cm]\nonumber
& +w_x[u_{xx}+v_{xy}+u_{xx}v_y+u_xv_{xy}-u_{xy}v_x-u_yv_{xx}]\\[.22cm]\nonumber
\kappa_{22}  =&~ u_y[v_{xy}w_y+v_xw_{yy}-v_{yy}w_x-(1+v_y)w_{xy}]\\[.1cm]
&+(1+v_y)[u_{yy}w_x+u_yw_{xy}-u_{xy}w_y-(1+u_x)w_{yy}]\\[.1cm]\nonumber
&+w_y[(1+v_y)u_{xy}+(1+u_x)v_{yy}-u_{yy}v_x-u_yv_{xy}]  \\[.22cm]\label{curve2}
\end{align}
\begin{align}\nonumber
\kappa_{12} =& ~  (1+u_x)[v_{xy}w_y+v_xw_{yy}-v_{yy}w_x-(1+v_y)w_{xy}] \\[.1cm]\nonumber
&+(1+v_y)[u_{xy}w_x+u_yw_{xx}-u_{xx}w_y-(1+u_x)w_{xy}]  \\[.1cm]\nonumber
&+u_y[v_{xx}w_y+v_xw_{xy}-v_{xy}w_x-(1+v_y)w_{xx}] \\[.1cm]\nonumber
&+v_x[u_{yy}w_x+u_yw_{xy}-u_{xy}w_y-(1+u_x)w_{yy}] \\[.1cm]\nonumber
&+w_x[(1+v_y)u_{xy}+(1+u_x)v_{yy}-u_{yy}v_x-u_yv_{xy}]  \\[.1cm]
&+w_y[u_{xx}+v_{xy}+u_{xx}v_y+u_xv_{xy}-u_{xy}v_x-u_yv_{xx}].\label{curve3}
\end{align}

In what follows we will use these expressions along with the inextensibility conditions in order to produce the potential energy from which we will obtain the equations of motion.
Let us first truncate these relations by dropping principal terms (second spatial derivatives) with coefficients of order higher than $\eta^2$ (according to the scaling in Assumption \ref{scaling}). After various cancellations this yields :
\begin{align}\label{curvaturenotsimplified}
\kappa_{11} =&~w_yv_{xx}+w_xu_{xx}-(1+u_x+v_y)w_{xx}\\[.1cm]
\kappa_{22}  =&~ u_{yy}w_x+w_yu_{yy}-(1+u_x+v_y)w_{yy}\\[.1cm] \label{cns*}
\kappa_{12} 
 =&~2[v_{xy}w_y+u_{xy}w_x-(1+u_x+v_y)w_{xy}].  
\end{align}

\begin{remark} Note that, in higher order models, we may wish to retain more terms in the above.\end{remark}
\subsection{Plate Model I: Three $\eta^2$-Inextensibility Conditions}
We begin with the potential energy expression from Section \ref{fullpotential}. We invoke the three full inextensibility conditions in \eqref{fullinext1}--\eqref{fullinext3} to eliminate $\varepsilon_{11}, \varepsilon_{22}, \varepsilon_{12}$ from the energetic expression. Also, as previously mentioned, we drop all $z^2\mu_{ij}$ terms in the expressions for $\varepsilon_{ij}$.
 This yields:
\begin{equation}
\label{potentialwithcurvature}
E_P =\frac{1}{2} \left[\frac{1}{12} \frac{E h^2}{(1-\nu^2)}\right] \int_{\Omega} \left[ \kappa^2_{11} + \kappa_{22}^2 + 2 \nu \kappa_{11} \kappa_{22} + \frac{1-\nu}{2}\kappa^2_{12} \right] d \Omega.
\end{equation}
We then invoke the curvature expressions in \eqref{curvaturenotsimplified}--\eqref{cns*} and input them in \eqref{potentialwithcurvature}. 

After truncating terms, we will obtain a potential energy (shown in the next section) for which we can enforce three {\em effective inextensibility relations} in \eqref{quadcon1}--\eqref{quadcon3}. Then, utilizing Hamilton's principle and following the same steps as in Section \ref{Beams}, we will obtain a straightforward PDE model with clear equations of motion. On the other hand, enforcing the three effective inextensibilty constraints results in three Lagrange multiplier variables that cannot be simultaneously eliminated in the full description of the system.

\subsubsection{Simplified Energies}\label{energiesplate1}

We employ the effective inextensibility constraints \eqref{quadcon1}--\eqref{quadcon3} in the curvature expressions given in \eqref{curvaturenotsimplified}. In particular, we can differentiate \eqref{quadcon1}--\eqref{quadcon3} variously, simplify, and rewrite $\kappa_{ij}$ solely in $w$

\begin{align*}
\kappa_{11} & = - w_{xx} \left[ 1 + \frac{1}{2} w_{x}^2 + \frac{1}{2} w_{y}^2 \right] \\[1em]
\kappa_{22} & = - w_{yy} \left[ 1 + \frac{1}{2} w_{x}^2 + \frac{1}{2} w_{y}^2 \right] \\[1em]
\kappa_{12} & =  - 2w_{xy} \left[ 1 +  \frac{1}{2} w_{x}^2 +  \frac{1}{2} w_{y}^2  \right].  \nonumber  
\end{align*}
We then form the appropriate products of $\kappa_{ij}$ as they appear in \eqref{potentialwithcurvature}. In line with the $\eta^2$ analysis of the beam, we retain in  $E_P$ only terms with coefficients up to and including $\eta^2$. (And thus we drop expressions of the form $[\partial_{x_i}w]^2[\partial_{x_j}w]^2\partial^2_{x_k}w$.) The result is the potential energy we will employ for the two principal models in this treatment:
\begin{equation}\label{goodPo}
E_P = \frac{D}{2} \int_0^{L_y} \int_0^{L_x} \left[ 1 + w_x^2 + w_y^2 \right] \left[ w_{xx}^2 + w_{yy}^2 + 2 \nu w_{xx}w_{yy} + 2(1- \nu) w^2_{xy} \right]  dx dy,
\end{equation}
where we have now denoted $D = \dfrac{1}{12} \dfrac{E h^2}{(1-\nu^2)}$.  We utilize the standard expression for the plate's kinetic energy (again with normalized mass density):
\begin{equation*}
E_K=\frac{1}{2}\int_0^{L_y}\int_0^{L_x}\left[u_t^2+v_t^2+w_t^2\right]dxdy.
\end{equation*}
\subsubsection{Equations of Motion}\label{eqmot1}
We introduce Lagrange multipliers $\lambda_i,~i=1,2,3$ acting to enforce the effective inextensibility constraints.
Let $\lambda_1$ be associated to the axial (effective) inextensibility condition \eqref{quadcon1}, $\lambda_2$ to the chord-wise condition \eqref{quadcon2}, and $\lambda_3$ to the shear constraint \eqref{quadcon3}. We invoke Hamilton's principle for the Lagrangian, written in the $\lambda_i$, $E_P$, and $E_K$, with virtual displacements $\delta u, \delta v,$ and $\delta w$. The calculations are involved but analogous to those in Section \ref{Beams} for the 1D dynamics. Essential boundary conditions are enforced for $u,v,$ and $w$ (and their virtual changes) on $\Gamma_W$.

In this presentation, since the variables $\lambda_i$ serve to enforce a constraint in the equations of motion, we will retain the equations \eqref{quadcon1}--\eqref{quadcon3} as part of the system:
\begin{equation} u_x+\frac{1}{2}w_x^2=0;~~~~v_y+\frac{1}{2}w_y^2=0;~~~~u_y+v_x+w_xw_y=0.\end{equation}
For the dynamics, we obtain:
\begin{align}
&u_{tt} + \partial_x \left( \lambda_1 \right) +  \partial_y \left( \lambda_3 \right) =0& \label{uttlambda}\\
& v_{tt} + \partial_y \left( \lambda_2 \right) +  \partial_x \left( \lambda_3 \right) =0 & \label{vttlambda}\\
&w_{tt} +{D}  \Delta[(1+|\nabla w|^2)\Delta w]-D\nabla\cdot [|\Delta w|^2\nabla w]  \nonumber \\
&~~~~~~~~+ \partial_{x} \left( \lambda_1 w_{x} \right) + \partial_{y} \left( \lambda_2 w_{y} \right) + \partial_{x} \left( \lambda_3 w_{y} \right) + \partial_{y} \left( \lambda_3 w_{x} \right) = 0. \label{wttlambda}&
\end{align}
The above system would be supplemented with appropriate initial displacement $w_0=w(x,y;0)$ and velocity $w_1(x,y) = w_t(x,y;0)$, from which the initial conditions for $u$ and $v$ can be inferred through the relationships \eqref{quadcon1}--\eqref{quadcon2}.
\begin{remark}
Letting $\bu = \langle u,v \rangle$ and $\ds \Lambda = \begin{bmatrix} \lambda_1 & \lambda_3 \\ \lambda_3 & \lambda_2 \end{bmatrix}$, we have a nice vectorial description of the system that makes for an apt comparison against other nonlinear plate equations (e.g., the full von Karman equations of motion \cite{ciarlet,lagnese,koch}):
\begin{align}
&\bu_{tt} + \text{div} \Lambda=0& \\
&w_{tt} + D\Big[ \Delta\big[(1+|\nabla w|^2)\Delta w\big]-\nabla\cdot \big(|\Delta w|^2\nabla w\big) \Big]+
\text{div}\Big(\Lambda \nabla w\Big)=0,&
\end{align} 
where the action of $\Lambda$ on the gradient is that of matrix multiplication in this presentation. 
\end{remark}

\subsubsection{Boundary Conditions}\label{BC1}
Here, we provide the boundary conditions for $w$ and for the $\lambda_i$. From these, the boundary conditions for $u$ and $v$ can be inferred. 
On the clamped edge $\Gamma_W$, we  again have the essential boundary condition
$$w=0;~~w_x= 0~~~ \text{on}~~ \Gamma_{W}.$$

The minimization of the Lagrangian and the arbitrariness of virtual displacements yield the {\em natural} boundary conditions from the potential energy. For the second order conditions in $w$, we obtain those of the standard linear free plate:
\begin{align*}
w_{xx} + \nu w_{yy} = 0~~~& \text{on}~~ \Gamma_{E}  \\
w_{yy} + \nu w_{xx} = 0~~~& \text{on}~~ \Gamma_{S} \\
w_{yy} + \nu w_{xx} = 0~~~& \text{on}~~ \Gamma_{N}. 
\end{align*}
The potential energy produces nonlinear forces (cubic-type)---and thus boundary conditions---along the free edges:
\begin{align*}
(1 - \nu) \left[ (1 + \nu)w_{x}w^2_{yy} - 2w_{x}w^2_{xy} - 4w_{y}w_{yy}w_{xy} \right]  - \left[ 1+w_x^2 + w_y^2 \right ] \left[ w_{xxx} + (2 - \nu) w_{yyx} \right] = 0~~~& \text{on}~~ \Gamma_{E} \\
(1 - \nu) \left[ (1 + \nu)w_{y}w^2_{xx} - 2w_{y}w^2_{xy} - 4w_{x}w_{xx}w_{xy} \right]  - \left[ 1+w_x^2 + w_y^2 \right ] \left[ w_{yyy} + (2 - \nu) w_{xxy} \right]= 0 ~~~ & \text{on}~~ \Gamma_{S}\\
(1 - \nu) \left[ (1 + \nu)w_{y}w^2_{xx} - 2w_{y}w^2_{xy} - 4w_{x}w_{xx}w_{xy} \right]  - \left[ 1+w_x^2 + w_y^2 \right ] \left[ w_{yyy} + (2 - \nu) w_{xxy} \right]= 0~~~ & \text{on}~~ \Gamma_{N}.
\end{align*}
\begin{remark} The above conditions are {\bf clearly not} the linear boundary conditions associated to the free plate, as presented in Remark \ref{linBC}.\end{remark}
The Lagrange multipliers $\lambda_i$, as variables, also have boundary conditions. These are readily obtained from Hamilton's principle along the free edges:
\begin{align*}
\lambda_1(x,y)=0~~\text{and}~~\lambda_3(x,y)=0~~~& \text{on}~~ \Gamma_{E} \\
\lambda_2(x,y)=0~~\text{and}~~\lambda_3(x,y)=0~~~& \text{on}~~ \Gamma_{S} \\
\lambda_2(x,y)=0~~\text{and}~~\lambda_3(x,y)=0~~~& \text{on}~~ \Gamma_{N}.
\end{align*}

Finally, we can express the boundary conditions for the Lagrange multipliers on the clamped edge $\Gamma_W$ using \eqref{uttlambda}--\eqref{vttlambda} and the equations above: 
\begin{align*}\small
&\lambda_1(0,y)  = \int_0^{L_x}\left[  u_{tt} +  \partial_y \left( \lambda_3 \right)\right]dx;~~
\lambda_2(0,y)  = \int_y^{L_y}\left[  v_{tt} +  \partial_x \left( \lambda_3 \right)\right]dy \bigg|_{x=0};~~
\lambda_3(0,y)  = \int_0^{L_x}\left[  v_{tt} +  \partial_y \left( \lambda_2 \right)\right]dx,&
\end{align*}
as well as 
\begin{align*}
\lambda_1(x,y) & = \int_x^{L_x}\left[  u_{tt} +  \partial_y \left( \lambda_3 \right)\right]dx\bigg|_{y=0}~~~ \text{on}~~ \Gamma_{S}, & ~~
\lambda_1(x,y) & = \int_x^{L_x}\left[  u_{tt} +  \partial_y \left( \lambda_3 \right)\right]dx\bigg|_{y=L_y} ~~~\text{on}~~ \Gamma_{N},\\[1em]
\lambda_2(x,y) & = \int_y^{L_y}\left[  v_{tt} +  \partial_x \left( \lambda_3 \right)\right]dy \bigg|_{x=L_x}~~~\text{on}~~ \Gamma_{E}.
\end{align*}

\subsubsection{Reduction of System \eqref{uttlambda}--\eqref{wttlambda}}
To be consistent with the beam system analysis, one may attempt to eliminate the $\lambda_i$, as well as the in-plane variables $u,v$ in the system. This would yield a dynamic equation in the principal elastic displacement $w$ only (as may be done for beam dynamics, and also as is possible in the case of von Karman's equations via the {\em Airy Stress Function} \cite{lagnese,ciarlet}). Yet, for the plate dynamics above, there seems to be no clear way to accomplish this. We opt, here, to eliminate $\lambda_1$ and $\lambda_2$, and to use \eqref{quadcon1} and \eqref{quadcon2} to write $u, v$ in terms of $w$. {\em This critically exploits the simplified nature of the span and chord-wise effective inextensibility constraints at the $\eta^2$-order.} In addition, we can see that, since the third (shear) constraint \eqref{quadcon3} does not permit us to explicitly solve for a displacement quantity, we will retain both $w$ and $\lambda_3$ in our reduced system.

To that end, we integrate \eqref{uttlambda} from $x$ to $L_x$ and utilize the boundary condition $\lambda_1(L_x, y)=0$. This yields:
\begin{equation}\label{lambda1expr}
\lambda_1(x,y) = \int_x^{L_x} \left[ u_{tt}  +  \partial_y \left( \lambda_3 \right) \right ] dx.
\end{equation}
Similarly, we integrate \eqref{vttlambda} from $y$ to $L_y$ and use the condition $\lambda_2(x,L_y)=0$. This gives:
\begin{equation} \label{lambda2expr}
\lambda_2(x,y) = \int_y^{L_y} \left[  v_{tt} +  \partial_x \left( \lambda_3 \right) \right] dy.
\end{equation}

Substituting \eqref{lambda1expr} and \eqref{lambda2expr} into \eqref{wttlambda} we obtain:
\begin{align}\label{plate11}
w_{tt} - D\Big [ \Delta[(1+|\nabla w|^2)\Delta w]-\nabla\cdot (|\Delta w|^2\nabla w) \Big] + \partial_{x} \left(w_{x} \int_x^{L_x} u_{tt}  \right) + \partial_{y} \left( w_{y} \int_y^{L_y} v_{tt}  \right)\nonumber \\[1em]+ \partial_{x} \left(w_{x} \int_x^{L_x} \partial_y \left(\lambda_3 \right)  \right) + \partial_{y} \left(w_{y} \int_y^{L_y} \partial_x \left(\lambda_3 \right)  \right) + \partial_{x} \left( \lambda_3 w_{y} \right) + \partial_{y} \left( \lambda_3 w_{x} \right) =0.
\end{align}
In-plane inertial expressions, $u_{tt}$ and $v_{tt}$, can be obtained, as in the case of the beam, by solving in \eqref{quadcon1}, \eqref{quadcon2}, and \eqref{quadcon3} and formally differentiating in time. Note that, using the essential boundary conditions at $x=0$ for each of $u,v,w$, we have:
\begin{align}
u(x,y)=&~-\frac{1}{2}\int_0^xw_x(\xi_1,y)^2d\xi_1\\
v(x,y)-v(x,0)=&~-\frac{1}{2}\int_0^yw_y(x,\xi_2)^2d\xi_2\\
v(x,y)=&~-\int_0^xw_x(\xi_1,y)w_y(\xi_1,y)d\xi_1-\int_0^{x}u_y(\xi_1,y)d\xi_1,
\end{align}
where we have suppressed the dependence on $t$ above. From which we have an expression for $v(x,y)$, and thence an expression for $u(x,y)$, where both depend only on the transverse variable $w$:
\begin{equation}
v(x,y) = \int_0^x\int_0^{\xi_1}w_x(\zeta,y)w_{xy}(\zeta,y) d\zeta d\xi_1-\int_0^xw_x(\xi_1,y)w_y(\xi_1,y)d\xi_1.
\end{equation}
From these we obtain the inertial expressions
\begin{align}\label{plate12}
u_{tt} = & ~-\int_0^x[w_{xt}^2+w_xw_{xtt}]d\xi_1 \\
v_{tt} = &~-\int_0^y[w_{yt}^2+w_yw_{ytt}]d\xi_2+v_{tt}(y=0). \label{plate13}
\end{align}

The above formulae in \eqref{plate11}--\eqref{plate13} showcase of the principal strengths of the $\eta^2$-order effective inextensibility constraints: namely, we obtain straight-forward equations of motion in two unknowns, with a clear connection to the $\eta^2$-order {\em beam dynamics} in Section \ref{beamsec1}.
On the other hand, owing to the structure of the third effective inextensibility constraint, \eqref{quadcon3}, it does not seem that systematic integrations and/or differentiations---coupled with the given boundary conditions---will allow a clean elimination of the $\lambda_3$ variable. Thus, any analytical treatment of this system would seem to require both the $w$ and $\lambda_3$ variables directly. 

\subsection{Model II: Partial Use of Effective Shear Constraint}
In this model, we retain the treatment of the potential energy in  Section \ref{energiesplate1}. To wit: we retain all three full inextensibility constraints \eqref{fullinext1}--\eqref{fullinext3} in eliminating the in-plane strains $\varepsilon_{ij}$ from the potential energy in \eqref{potentialfoundamental}. After this, we truncate all three inextensibility constraints to quadratic order (as in \eqref{quadcon1}--\eqref{quadcon3}). Finally, we truncate the expanded potential energy as before to arrive at the potential energy expression in \eqref{goodPo}. On the other hand, in developing the equations of motion, we  elect (with foresight coming from those issues associated to $\lambda_3$ above) to {\em enforce only the first two $\eta^2$-effective inextensibility conditions}, \eqref{quadcon1} and \eqref{quadcon2} through Lagrange multipliers $\lambda_1, \lambda_2$. Thus, in this model, shear inextensibility is only implicitly enforced from the point of view of the choice of the potential energy. We utilize only two Lagrange multipliers in the derivation of the equations of motion.

In this presentation, the constraint variables $\lambda_1,~\lambda_2$ serve to enforce the span and chord-wise quadratic effective inextensibility constraints.  Accordingly, we retain these as part of the system:
\begin{align}\label{refinextnow}
u_{x} + \frac{1}{2} w_{x}^2 =0;~~~~v_{y} + \frac{1}{2} w_{y}^2 =0.
\end{align}
For the dynamics, Hamilton's principle yields identical nonlinear terms in $w$, owing the use of $E_P$ as in Section \ref{eqmot1}. Moreover, we do not retain any reference to the shear constraint in the equations, as there is no $\lambda_3$ present. The unforced equations of motion are then:
\begin{align}\label{reduced1}
&u_{tt} + \partial_x \left( \lambda_1 \right)    =0&\\\label{reduced2}
& v_{tt} + \partial_y \left( \lambda_2 \right)  =0 &\\ \label{reduced3}
&w_{tt} + D \left [ \Delta[(1+|\nabla w|^2)\Delta w]-\nabla\cdot (|\Delta w|^2\nabla w) \right] + \partial_{x} \left( \lambda_1 w_{x} \right) + \partial_{y} \left( \lambda_2 w_{y} \right)  = 0.&
\end{align}
As before, the equations would be supplemented with appropriate initial displacement $w_0=w(x,y;0)$ and velocity $w_1(x,y) = w_t(x,y;0)$, from which the initial conditions for $u$ and $v$ can be inferred. 

\subsubsection{Boundary Conditions}
The boundary conditions for $w$ (which then yield conditions for $u$ and $v$ through \eqref{refinextnow}) are identical to Section \ref{BC1}. 
The boundary conditions of the Lagrange multipliers $\lambda_1$ in this case are:
\begin{align*}
\lambda_1(x,y) & = \int_0^{L_y}u_{tt} dy~~~\text{on }\Gamma_{W},&&\lambda_1(x,y)  =0~~~\text{on}~~ \Gamma_{E},\\[1em]
~\lambda_1(x,y) & = \int_x^{L_x} u_{tt}dx\bigg|_{y=0}~~~ \text{on}~~ \Gamma_{S},
&&\lambda_1(x,y) = \int_x^{L_x} u_{tt} dx\bigg|_{y=L_y} ~~~\text{on}~~ \Gamma_{N}.
\end{align*}

And for $\lambda_2$:
\begin{align*}
\lambda_2(x,y) & = \int_y^{L_y} v_{tt}dy \bigg|_{x=0}~~~ \text{on}~~ \Gamma_{W}, &&\lambda_2(x,y)  = \int_y^{L_y} v_{tt}dy \bigg|_{x=L_x}~~~\text{on}~~ \Gamma_{E}, \\[1em] \lambda_2(x,y) & =0~~~ \text{on}~~ \Gamma_{S}, && \lambda_2(x,y)=0~~~ \text{on}~~ \Gamma_{N}.
\end{align*}

\subsubsection{Reduction of System \eqref{reduced1}--\eqref{reduced3}}
One of the primary benefits for this system, discussed further in Section \ref{Discussion}, is that the Lagrange multiplier variables can be fully eliminated (as with the inextensible beam). Indeed, one obtains from \eqref{reduced1} and \eqref{reduced2} that
$$\lambda_1(L_x,y)-\lambda_1(x,y)=-\int_x^{L_x}u_{tt}d\xi;~~~~\lambda_2(x,y)-\lambda_2(x,0)=-\int_0^yv_{tt}d\zeta.$$
At this point we can invoke the boundary conditions for $\lambda_i$, as seen above on the appropriate edge of the plate, to conclude
\begin{equation}
\lambda_1(x,y)=\int_x^{L_x}u_{tt}d\xi;~~~~\lambda_2(x,y)=-\int_0^yv_{tt}d\zeta.
\end{equation}

These quantities may be substituted directly into the equations of motion for $w$, which results in the following {\em closed} system---with no reference to the Lagrange variables:
\begin{align}\label{citing1}
&w_{tt} + D\Delta[(1+|\nabla w|^2)\Delta w]-D\nabla\cdot [|\Delta w|^2\nabla w] + \partial_{x} \big( w_{x} \int_x^{L_x}u_{tt}d\xi  \big) - \partial_{y} \big( w_{y}\int_0^yv_{tt}d\zeta \big)  = 0&\\ 
&\hskip2cm u_{x} + \frac{1}{2} w_{x}^2 =0;~~~v_{y} + \frac{1}{2} w_{y}^2 =0.\label{citing2}
\end{align}
A complete description, then, would again provide the relevant initial and boundary conditions from the principal variable here, $w$. 

\subsection{Model III:  Complete Omission of The Shear Constraint}
In natural succession from the previous models, one may inquire: \begin{quote} What happens if we omit the third (shear) constraint in the derivation?\end{quote} This is a reasonable subsequent step, as we have in the previous model only partially made use of the shear constraint (both in its full form, as well as in its simplified, quadratic effective form). For the model in this section, we refrain entirely from making mention to a shear constraint. With only span and chord-wise inextensibility enforced, we obtain the equations of motion. In particular, this showcases how convoluted the equations of motion become and demonstrates the usefulness of the shear constraint.

\subsubsection{Addressing The Potential Energy}
The two constraints approach forbids elongation in the $x$ and $y$ axes but {\em permits} shear strain. This translates into $ {\varepsilon}_{11} = {\varepsilon}_{22}=0$, but we will {\em not} take $ \varepsilon_{12} = 0$. Immediately it is clear we will have more terms in the equations of motion. As before, after invoking the full inextensibility conditions $\varepsilon_{11}=0$ and $\varepsilon_{22}=0$, we will compute the associated potential energy and then truncate to a particular order. 

Applying only two inextensibility constraints to eliminate $\varepsilon_{11}, \varepsilon_{22}$  into the bulk strain expressions \eqref{bulk1}--\eqref{bulk3}  yields:
\begin{align*}
\h \varepsilon_{11} & = -z \left[w_{xx}  \left(1+ u_{x} + v_{y} \right) -  w_{x}  u_{xx} - w_{y} v_{xx} \right] \nonumber \\[1em]
 \h \varepsilon_{22} & = -z \left[   w_{yy}  \left(1+ u_{x} + v_{y} \right) - w_{x}  u_{yy} - w_{y} v_{yy} \right] \\[1em]
 \h \varepsilon_{12} & = u_{y} +v_{x} +u_{x}u_{y} +v_{x}v_{y} +w_{x}w_{y} - 2z \left[w_{xy}  \left(1+ u_{x} + v_{y} \right) - w_{x} u_{xy} - w_{y}v_{xy} \right]. 
\end{align*}
We then invoke the $\eta^2$-order hypothesis to obtain the curvature expressions \eqref{curvaturenotsimplified}--\eqref{cns*}.  We must also truncate $${\varepsilon}_{12}=u_{y} +v_{x} +u_{x}u_{y} +v_{x}v_{y} +w_{x}w_{y} \approx u_{y} +v_{x} +w_{x}w_{y},$$ under the same $\eta^2$-hypothesis. And, as before, we invoke the (two) $\eta^2$-effective inextensibility constraints 
$$u_x=-\frac{1}{2}w_x^2,~~~v_y=-\frac{1}{2}w_y^2,$$
from which we can differentiate to solve for $u_{xx}$ and $v_{yy}$. 
\begin{remark} Without creating nonlocal conditions---and invoking boundary conditions as yet undetermined---we can no longer obtain expressions for $v_x, u_y$ as before. We note that this is a key point of distinction when we have imposed no notion of shear inextensibility. \end{remark}
Upon simplifying the terms as described above, we obtain:
\begin{align*}
\hat \varepsilon_{11} & = z \left[ - w_{xx} \left( 1 + \frac{1}{2} w_{x}^2 - \frac{1}{2} w_{y}^2 \right) + w_{y}v_{xx} \right] \nonumber \\[1em]
\hat \varepsilon_{22} & =  z \left[ - w_{yy} \left( 1 - \frac{1}{2} w_{x}^2 + \frac{1}{2} w_{y}^2 \right) + w_{x}u_{yy} \right]  \\[1em]
 \hat \varepsilon_{12} & =  u_{y} + v_{x} + w_{x}w_{y} - 2z w_{xy} \left[  1 +  \frac{1}{2} w_{x}^2 +  \frac{1}{2} w_{y}^2  \right]. \nonumber  
\end{align*}
Then, integrating in $z$ and simplifying, with coefficients:
\begin{align*}
\int_{-h/2}^{h/2}[\hat{\varepsilon}_{11}^2+\hat{\varepsilon}_{22}^2]dz= & ~\dfrac{h^3}{12}\Big[w_y^2v_{xx}^2+w_x^2u_{yy}^2-2w_yv_{xx}w_{xx}-2w_xu_{yy}w_{yy} \\
&~~~~~ -w_y(w_x^2-w_y^2)v_{xx}w_{xx}+w_x(w_x^2-w_y^2)u_{yy}w_{yy} \\
&~~~~~ + w_{xx}^2\Big(1+(w_x^2-w_y^2)+\frac{1}{4}(w_x^2-w_y^2)^2\Big)+w_{yy}^2\Big(1-(w_x^2-w_y^2)-\frac{1}{4}(w_x^2-w_y^2)^2\Big)\Big]\\
2\nu \int_{-h/2}^{h/2}[\hat{\varepsilon}_{11}\hat{\varepsilon}_{22}]dz= & ~ \frac{\nu h^3}{6}\Big[v_{xx}u_{yy}w_xw_y+w_{xx}w_{yy} w_xw_{xx}u_{yy}-w_yv_{xx}w_{yy}\\
&~~~~~~~~~  -\frac{1}{4}w_{xx}w_{yy}\big(w_x^2-w_y^2)^2 + \frac{1}{2}\big(w_x^2-w_y^2)(w_yv_{xx}w_{yy}-w_xw_{xx}u_{yy}\big)\Big]\\
\dfrac{1-\nu}{2} \int_{-h/2}^{h/2}[\hat{\varepsilon}_{12}^2]dz= & ~ \frac{h(1-\nu)}{2}(u_y+v_x+w_xw_y)^2+\frac{h^3(1-\nu)}{6}\Big(1+\frac{1}{2}w_x^2+\frac{1}{2}w_y^2\Big)^2w_{xy}^2.
\end{align*}

Recalling the full potential energy expression \eqref{potentialfoundamental}, we input the above expressions for the strains. Since we are operating at the $\eta^2$-level here for coefficients, we discard {\em any terms} with coefficients scaled by $\eta^3$ or higher. Note that, in the (new) situation where first derivative terms appear independent of second derivative terms, we choose to retain. This yields {\small
\begin{align*}
E_P = &~ \frac{1}{2}\frac{E}{1-\nu^2}\int_{\Omega}
\dfrac{h^3}{12}\Big[w_y^2v_{xx}^2+w_x^2u_{yy}^2-2w_yv_{xx}w_{xx}-2w_xu_{yy}w_{yy}+ w_{xx}^2\Big(1+w_x^2-w_y^2\Big)+w_{yy}^2\Big(1-w_x^2+w_y^2\Big)\Big]\\
&+\frac{\nu h^3}{6}\Big[v_{xx}u_{yy}w_xw_y+w_{xx}w_{yy} w_xw_{xx}u_{yy}-w_yv_{xx}w_{yy}\Big]+\frac{h^3(1-\nu)}{6}\Big(1+w_x^2+w_y^2\Big)w_{xy}^2\\
&+ \frac{h(1-\nu)}{2}(u_y+v_x+w_xw_y)^2 ~d\Omega.
\end{align*}}
\begin{remark}The terms in final line above represent an interesting contribution, both in terms of their ``$h$" scaling, i.e., thickness, as well as being detached from principal second derivative terms. \end{remark}
\begin{remark} For reference, the terms we discarded above include: {\small
\begin{align*}
&\dfrac{h^3}{12}\Big[w_x(w_x^2-w_y^2)u_{yy}w_{yy} -w_y(w_x^2-w_y^2)v_{xx}w_{xx}+\frac{w_{xx}^2}{4}(w_x^2-w_y^2)^2-\frac{w_{yy}^2}{4}(w_x^2-w_y^2)^2\Big ]&\\
&-\dfrac{\nu h^3}{24}w_{xx}w_{yy}\big(w_x^2-w_y^2)^2+ \frac{\nu h^3}{12}\big(w_x^2-w_y^2)(w_yv_{xx}w_{yy}-w_xw_{xx}u_{yy}\big)+\frac{h^3(1-\nu)}{6}\Big(\frac{1}{2}w_x^2w_y^2+\frac{1}{4}w_x^4+\frac{1}{4}w_y^4\Big)w_{xy}^2.&
\end{align*}} \end{remark}

Now, invoking the definition of the constant $D$, the potential energy can then be written:
\begin{align}
\label{potentialfoundamental2final}\nonumber
E_P =\frac{6D}{h} \int_{0}^{L_y} \int_0^{L_x} \bigg \{ \frac{h^2}{12} & \big [ w^2_{xx}+w^2_{yy}+2\nu w_{xx}w_{yy}+ 2(1-\nu) w^2_{xy} \left(  1 + w_{x}^2 +   w_{y}^2  \right)  + (w_x^2-w_y^2)(w_{xx}^2-w_{yy}^2)  \\&  - 2 w_{y} w_{xx} v_{xx}  -2 w_{x} w_{yy} u_{yy}- 2 \nu \left(   w_{y} w_{yy} v_{xx} +w_{x}w_{xx}u_{yy} \right)  \big ]  \nonumber
\\& + \frac{1-\nu}{2} \big[u_y+v_x+w_xw_y]^2 \bigg \}  dx dy.
\end{align}

\subsubsection{Equations of Motion}
As before, $\lambda_1$ is used to enforce the span-wise effective constraint \eqref{quadcon1}  and $\lambda_2$ the chord-wise \eqref{quadcon2}. Utilizing the arbitrariness of the relevant virtual changes, we can gather the equations of motion. We recover both effective constraints \eqref{quadcon1} and \eqref{quadcon2} via the associated Lagrange multipliers $\lambda_1$ and $\lambda_2$ in Hamilton's principle. 
Following the procedure of the previous sections, the equations of motion for the in-plane displacements: {
\begin{align*}
u_{tt} + \partial_x \left( \lambda_1 \right) -D \left [ w_{xyy}w_{yy} +2w_{xy}w_{yyy}  +w_{x} \partial^4_y w + \nu w_{xyy} w_{xx} + 2\nu w_{xy} w_{xxy} + \nu w_{x} w_{xxyy} \right ] \\
 - \frac{12D}{h^2}(1- \nu)  \left [u_{yy} + w_{x}w_{yy} \right ] =0
\end{align*}
\begin{align*}
v_{tt} + \partial_y \left( \lambda_2 \right) - 2D\left [w_{yxx}w_{xx} +2w_{yx}w_{xxx}  + w_{y} \partial^4_x w + \nu w_{yxx} w_{yy} + 2\nu w_{yx} w_{yyx} + \nu w_{y} w_{yyxx} \right ] \\
-\frac{12D}{h^2}(1- \nu) \left [v_{xx} + w_{xx}w_{y}  \right ] =0.
\end{align*}}

We rewrite these, and include the equation for $w$:
\begin{align}
u_{tt} + \partial_x \left( \lambda_1 \right)- \frac{12D}{h^2}(1- \nu) \left [u_{yy} + w_{x}w_{yy} \right ] - 2D \partial_y^2\left [w_x( w_{yy}+\nu w_{xx})\right ] 
 =0\\
v_{tt} + \partial_y \left( \lambda_2 \right) - \frac{12D}{h^2}(1- \nu) \left [v_{xx} + w_{xx}w_{y}  \right ]   -2D \partial_x^2\left[w_y(w_{xx}+\nu w_{yy}) \right ] =0
\end{align}
\small{
\begin{align*}
&w_{tt} + \partial_{x} \left( \lambda_1 w_{x} \right) + \partial_{y} \left( \lambda_2 w_{y} \right) -D \Big [ \partial^4_x w (1 + w_x^2 -w^2_y) + \partial^4_y w (1 - w_x^2 + w^2_y) +2w_{xxyy} (1 + w_x^2 + w^2_y) \\[1em]
& -  \nu w_{xxyy} (w_x^2 + w^2_y) -w_{x} \partial_y^4 u - w_{y} \partial_x^4 v + 4w_{x}w_{xx}w_{xxx} + 4w_{y}w_{yy}w_{yyy} - w_{x}w_{yy}w_{xyy} - w_{y}w_{xx}w_{yxx} \\[1em] 
& + (4-2\nu)w_{x}w_{xy}w_{xxy} + (4-2\nu)w_{y}w_{xy}w_{yyx} - 4w_{x}w_{xy}w_{yyy} - 4w_{y}w_{yx}w_{xxx} - 2w_{xy}v_{xxx} - 2w_{xy}u_{yyy}\\[1em] 
& + 4(1-\nu)w_{x}w_{xx}w_{xyy} + 4(1-\nu)w_{y}w_{yy}w_{xxy}  + w_{xx}^3 + w_{yy}^3 +w_{xx}w_{yy}^2 + w_{yy}w_{xx}^2  -(1 + 3\nu)w_{xx}w_{xy}^2 - (1 + 3\nu)w_{yy}w_{xy}^2 \Big ] \\[1em] 
& +\frac{6D}{h^2}  (1 - \nu) \Big [ w_{xx}w_{y}^2 + 2w_{x}w_{y}w_{xy} + 2u_{y}w_{xy} + v_{xx}w_{y} + 2v_{x}w_{xy} + w^2_{x}w_{yy} + u_{yy}w_{x} \Big ]  =0.
\end{align*}}

\subsubsection{Boundary Conditions}
On the clamped edge $\Gamma_W$, we have:
\begin{align*}
w=0;~~~w_x= 0;~~~u=0;~~~v=0. 
\end{align*}
For the second order conditions, we have:

On $\Gamma_{E}$:
\begin{align*}
&w_{xx} + \nu w_{yy}=0,~~
&w_{xx} (1 + w_x^2 -w^2_y) -w_{y}v_{xx} + \nu w_{yy} - \nu w_{x}u_{yy} =0.
\end{align*}

On $\Gamma_S$ and $\Gamma_N$:
\begin{align*}
&w_{yy} + \nu w_{xx} =0,~~
&w_{yy} (1 - w_x^2 + w^2_y) -w_{x}u_{yy}+ \nu w_{xx}  - \nu w_{y}v_{xx}   =0. 
\end{align*}
For the third order conditions we have:\\

On $\Gamma_E$:
\begin{align*}
& \frac{h^2}{6} w_{y} \left [ w_{xxx} + \nu  w_{yyx}\right  ]   + (1- \nu)\left [ v_{x} + w_{x}w_{y} + u_{y} \right ] =0 \\[1em]
& \frac{h^2}{6} \Big \{  -w_{x}w^2_{xx} - w_{x}w^2_{yy} -w_{yy}u_{yy}  - (2 - \nu)w_{x} w^2_{xy}  - w_{xxx} (1 + w_x^2 -w^2_y)  + 2w_{y}w_{xx}w_{yx} \\[1em]
&~~ + w_{yx}v_{xx} + w_{y}v_{xxx} - \nu w_{yyx}  - \nu w^2_{x}w_{xyy} -2(1 - \nu) \left[ w_{xyy} (1 + w_x^2 + w_y^2) + 2w_yw_{xy} w_{yy}   \right]\Big \} \\[1em]
&~~~~ +(1 - \nu) \left[w_{x}w_{y}^2 + u_{y}w_{y} + v_{x}w_{y} \right ]  =0.
\end{align*}

On $\Gamma_S$ and $\Gamma_N$:
\begin{align*}
&\frac{h^2}{6}w_{x}\left [ w_{yyy} + \nu w_{xxy} \right ] +(1- \nu)  \left [u_{y} + w_{x}w_{y} + v_{x} \right ]=0 \\[1em]
&\frac{h^2}{6} \Big \{ - w_{y}w^2_{yy} - w_{y}w^2_{xx}  -w_{xx}v_{xx} -  (2-\nu)w_{y} w^2_{xy} -w_{yyy} (1 - w_x^2 + w^2_y) + 2w_{x}w_{yy}w_{yx}\\[1em]
&~~ + w_{yx}u_{yy} + w_{x}u_{yyy} - \nu w_{xxy}   - \nu w^2_{y}w_{yyx}  -2(1 - \nu)\left[ w_{xxy} (1 + w_x^2 + w_y^2) +2w_{x}w_{xy} w_{xx} \right] \Big \}\\[1em]
&~~~~ + (1 - \nu) \left[w^2_{x}w_{y} + u_{y}w_{x} + v_{x}w_{x} \right ] =0. 
\end{align*}

The boundary conditions of the Lagrange multipliers $\lambda_i$ are as follow:
\begin{align*}
\lambda_1(x,y)=0~~~ \text{on}~~ \Gamma_{E};~~
\lambda_2(x,y)=0~~~ \text{on}~~ \Gamma_{S};~~
\lambda_2(x,y)=0~~~ \text{on}~~ \Gamma_{N};
\end{align*}
\begin{align*}
\lambda_1(x,y) & = \int_{0}^{L_y} \left \{ u_{tt} - \frac{12D}{h^2}(1- \nu) \left [u_{yy} + w_{x}w_{yy} \right ] - 2D \partial_y^2\left [w_x( w_{yy}+\nu w_{xx})\right ] \right \} dy \bigg |_{x=0}~~~ \text{on}~~ \Gamma_{W} \\[1em]
\lambda_1(x,y) & = \int_{x}^{L_x} \left \{ u_{tt} - \frac{12D}{h^2}(1- \nu) \left [u_{yy} + w_{x}w_{yy} \right ] -2D \partial_y^2\left [w_x( w_{yy}+\nu w_{xx})\right ] \right \} dx \bigg |_{y=0} ~~~ \text{on}~~ \Gamma_{S} \\[1em]
\lambda_1(x,y) & = \int_{x}^{L_x} \left \{ u_{tt} - \frac{12D}{h^2}(1- \nu) \left [u_{yy} + w_{x}w_{yy} \right ] - 2D \partial_y^2\left [w_x( w_{yy}+\nu w_{xx})\right ] \right \} dx \bigg |_{y=L_y} ~~~ \text{on}~~ \Gamma_{N} \\[1em]
\lambda_2(x,y) & = \int_{y}^{L_y}  \left \{ v_{tt} - \frac{12D}{h^2}(1- \nu) \left [v_{xx} + w_{xx}w_{y}  \right ]   - 2D \partial_x^2\left[w_y(w_{xx}+\nu w_{yy}) \right ] \right \} dy \bigg |_{x=0}~~~ \text{on}~~ \Gamma_{W} \\[1em]
\lambda_2(x,y) & = \int_{y}^{L_y}  \left \{ v_{tt} - \frac{12D}{h^2}(1- \nu) \left [v_{xx} + w_{xx}w_{y}  \right ]   - 2D \partial_x^2\left[w_y(w_{xx}+\nu w_{yy}) \right ] \right \} dy \bigg |_{x=L_x}~~~ \text{on}~~ \Gamma_{E} .
\end{align*}

\subsection{Higher Order Models}\label{higherorder}
We conclude the central part of the treatment of inextensible plates by addressing the natural question of including higher order terms in various modelling steps. We see that in the approaches above that, unlike for the beam, there are critical junctures where the order (of truncation) affects the inextensibility considerations, as well as (though not independent of) the analysis of the potential energy. On the one hand, it is possible to  address some of the shortcomings of the previous three models via the inclusion of higher order effects; on the other hand, as we have already seen, the simplest possible (and lowest order) truncations already yield models which are notably complex.

\subsubsection{Higher Order Inextensibility} 

We recall the ``full" inextensibility constraints: 
\begin{align}
(1+u_x)^2+v_x^2+w_x^2=&~1 \label{const1} \\ 
u_y^2+(1+v_y)^2+w_y^2=&~1 \label{const2} \\
u_y+v_x+u_xu_y+v_xv_y+w_xw_y=&~0.
\end{align}
We can, of course, expand the first two to read:
\begin{align}
2u_x+u_x^2+v_x^2+w_x^2=~0;&~~~~2v_y+v_y^2+u_y^2+w_y^2=~0.
\end{align}
\begin{remark} If we were to combine the span and chord constraints, we would have:
\begin{align}
\nabla \cdot \langle u, v \rangle=&~-\frac{1}{2}\big[|\nabla u|^2+|\nabla v|^2+|\nabla w|^2\big] \\ 
u_y+v_x=&~-[u_xu_y+v_xv_y+w_xw_y].
\end{align}
\end{remark}
Now, we may proceed as we did before in Assumption \ref{scaling} (and the discussion thereafter), and  retain higher order terms (up to $\eta^4$) in \eqref{taylor1}--\eqref{taylor2}. In this case, we would obtain:
\begin{align}
u_x=& ~-\frac{1}{2}[w_x^2+v_x^2]-\frac{1}{8}w_x^4 \\
v_y=& ~-\frac{1}{2}[w_y^2+u_y^2]-\frac{1}{8}w_y^4 \\
u_y+& v_x=-u_xu_y-v_xv_y-w_xw_y.
\end{align}
On the other hand, more closely following the logic in 1D, we would have from \eqref{const1} that $w_x^2 \sim \eta^2$, and also then that $u_x \sim \eta^2$, so we would discard {\bf only} $u_x^2$ , but not $u_y^2$. Similarly, from \eqref{const2}, we would discard $v_y^2$ only but not $v_x^2$. 
We can formalize this into a separate hypothesis.
\begin{assumption}\label{scaling2} Assume that $\partial_{x_i}w, \partial_yu, \partial_xv \sim \eta$, and assume that $[\partial_{x} u], [\partial_{y} v] \sim \eta^2$, for all $i=1,2$.\end{assumption} 
Formulating the reduced $\eta^4$ constraints from the Taylor expansions in \eqref{taylor1}---\eqref{taylor2}, we obtain:
\begin{align}
u_x=&  ~-\frac{1}{2}v_x^2-\frac{1}{2}w_x^2-\frac{1}{4}v_x^2w_x^2-\frac{1}{8}v_x^4-\frac{1}{8}w_x^4\\
v_y=&~-\frac{1}{2}u_y^2-\frac{1}{2}w_y^2-\frac{1}{4}u_y^2w_y^2-\frac{1}{8}u_y^4-\frac{1}{8}w_y^4\\
u_y+& v_x=-u_xu_y-v_xv_y-w_xw_y.
\end{align}

\subsubsection{Higher Order Potential Energy}
With either of the choices in the previous section, we might develop an $\eta^4$ potential energy. (This is principally a distinction between taking Assumption \ref{scaling} and Assumption \ref{scaling2}.) In doing so, we would revisit the analysis as before, beginning with the potential energy expression
\begin{equation}
E_P = \frac{1}{2} \frac{E}{1-\nu^2} \int_{\Omega}\int_{-h/2}^{h/2} \left[ \h \varepsilon^2_{11} + \h \varepsilon_{22}^2 + 2 \nu \h \varepsilon_{11} \h \varepsilon_{22} + \frac{1-\nu}{2}\h \varepsilon^2_{12} \right] dz d \Omega.
\end{equation}
We would reconsider the full strains, discarding the terms which are quadratic in $z$ (as in Assumption \ref{z}):
\begin{align}
 \h \varepsilon_{11} =  { \varepsilon} _{11} +z \kappa_{11};~~~
\h \varepsilon_{22}  =  { \varepsilon} _{22} +z \kappa_{22};~~ \h \varepsilon_{12} & =  { \varepsilon} _{12} +z \kappa_{12}. 
\end{align}
Recalling the strain-displacement relations
\begin{align}\label{strain1*}
\varepsilon_{11} & = u_{x} + \frac{1}{2} \left[ u_x^2  + v^2_x + w_x^2 \right]\\[1em]
\varepsilon_{22} & = v_{y} + \frac{1}{2} \left[ u_y^2  + v^2_y + w_y^2 \right]\\[1em]\label{strain3*}
\varepsilon_{12} & =u_y + v_x + u_x u_y + v_x v_y + w_x w_y,
\end{align}
we may then choose to enforce inextensibility by taking each ${\varepsilon}_{ij}=0, i,j=1,2$, or by only taking ${ \varepsilon}_{11}={ \varepsilon}_{22}=0$. In either case, we must then implement curvature expressions which are themselves ``accurate" up to the order of $\eta^4$. For the sake of space, we do not reproduce the curvature expressions from \eqref{curve1}--\eqref{curve3} here, but suffice to say that one may implement either Assumption \ref{scaling} or Assumption \ref{scaling2} {\em up to order $\eta^4$}. We demonstrate with \eqref{curve1}:
\begin{align} \label{curve1*}
\kappa_{11} =&~(1+u_x)[v_{xx}w_y+v_xw_{xy}-v_{xy}w_x-(1+v_y)w_{xx} ]\\[.1cm]\nonumber
&+v_x[u_{xy}w_x+u_yw_{xx}-u_{xx}w_y-(1+u_x)w_{xy}]\\[.1cm]\nonumber
& +w_x[u_{xx}+v_{xy}+u_{xx}v_y+u_xv_{xy}-u_{xy}v_x-u_yv_{xx}].
\end{align}
Recall that under Assumption \ref{scaling}, up to order $\eta^2$, we obtained: 
$$\kappa_{11} =~w_yv_{xx}+w_xu_{xx}-(1+u_x+v_y)w_{xx}.$$

Now, including all terms up to order $\eta^4$ (from either Assumption \ref{scaling} or Assumption \ref{scaling2}), we should retain all terms in \eqref{curve1*}. Indeed,  coefficients on second order terms are at or below $\eta^4$ order in their coefficients. Additionally, unlike the resulting analysis from the $\eta^2$ inextensibility assumption, we cannot explicitly solve for $u_{x_ix_j}$ and $v_{x_ix_j}$ to cleanly simplify the expressions in $\kappa_{11}$. Thus, pursuing an $\eta^4$-order model from either Assumption \ref{scaling} or Assumption \ref{scaling2} will require retaining expressions for $\kappa_{ij}$ in their entirety, resulting in a rather perilous expression for the potential energy $E_P$.

We do not pursue this line further here. We have included the discussion above to demonstrate a general method of building a model that consistently implements an order hypothesis across the inextensibility constraints. Moreover, we observe that it is possible to do so to ever-increasing order. On the other hand, as we can already see through Models 1--3 above, the degree of complexity in the PDE models, boundary conditions, and enforcement of the inextensibility constraints through Lagrange multipliers, is notable at the level of even $\eta^2$.

\section{Discussion and Future Work}\label{Discussion}

We now provide some  discussion of the plate models presented above to supplement the remarks contained in the sections. 
\begin{itemize}
\item The principal benefit of building $\eta^2$-plate models is they maintain $u$ and $v$ as decoupled variables, in that one can solve for them in terms of $w$. (This was the case of the beam analysis.) This, of course, is tremendously useful in developing a clear and tractable set of equations of motion for performing subsequent analyses.	
\item The first point of departure from the beam theory presented in Section \ref{Beams} is the emergence of nonlinear boundary conditions (Section \ref{BC1}); this is a necessary byproduct of invoking Hamilton's principle with the chosen potential energy (very similar to what occurs in \cite{htw} in an extensible situation). For any future theory of existence and uniqueness of solutions for the inextensible plate mirroring that of the beam \cite{DW}, nonlinear boundary conditions will be a central issue. This is especially true, as the nonlinear boundary conditions emerge in the notoriously troublesome higher-order plate conditions, and differ markedly from the linear theory. In any approach making using mode functions, one has to be cognizant of this point, in particular, from the point of view of convergence analysis and associated error near the free edges. 
	\item Furthermore, even the {\em linear mode functions} associated to 2D cantilevered plates are notoriously challenging \cite{beamforplates}, far beyond the simple cantilever mode functions in 1D \cite{hhww,htw}.  Moreover, with the aforementioned appearance of nonlinear boundary conditions, linear mode functions will not satisfy relevant inextensible plate boundary conditions, even if one could explicitly express such linear mode functions through separation \cite{Leissa,Blevins}.  It is common in the engineering community to assume mode shapes are multiples of the 1D beam mode shapes, although for plates with free edges this is not exact \cite{beamforplates}.
\item The second major point of departure is the inability to neatly and cleanly resolve the equations through the elimination of the Lagrange multipliers in Plate Models I and III. The benefit of Model II is that it yields equations of motion that are written explicitly in the transverse deflection $w$; on the other hand, it has an inherent discrepancy between the enforcement of the shear inextensibility constraint in the potential energy (taken) and the Lagrange analysis (ignored as a constraint).
\item The dissertation \cite{kevinthesis} demonstrated a numerical implementation of Model I described in this paper using a Rayleigh-Ritz global modal method. It was found that the model was less stiff than a nonlinear model analyzed by ANSYS, a commercial Finite Element Method solver. Model III was proposed in the dissertation, and the same method was used to implement that model numerically. Using the potential energy associated to Model III resulted in a plate which was far stiffer than the one modeled in ANSYS.
		\item In the discussion of cantilever large deflections, it is not physically appropriate to discard $u_{tt}$ or $v_{tt}$ in the modeling as one would do in the case of the scalar von Karman dynamics \cite{lagnese,springer}. Such a decision would, of course, simplify the model; however, as nonlinear inertial effects are paramount here, doing so would be a poor modeling choice. For a cantilevered beam it is easily observed that these ``acceleration terms" are on the order of nonlinear stiffness terms for the first mode \cite{dowell4,inext2}. In general, this will be true for the cantilevered plate, though this is not the case for beams or plates which are restrained on opposing boundaries (again, see the von Karman plate theory \cite{lagnese,ciarlet}).

\item As a first mathematical step, one should solve the stationary problem(s) associated to the potential energy $E_P$ in \eqref{goodPo}. This should be done from the point of view of internal and boundary loading. Some preliminary numerics---discussed above---have considered edge loading, but no thorough numerical or analytical studies have been undertaken. Mathematically, a first step will be developing a theory of strong solutions for the stationary version of Model I (dropping time derivatives and Lagrange multiplier terms). The resulting model is spatially-quasilinear with nonlinear boundary conditions.

\item It is apparent from the analysis in Model III and the discussion in Section \ref{higherorder} that higher order methods can and should be developed. However, owing to the apparent complexity of such models, their development should be largely tailored to the context of what is to be modeled---for instance, paying close attention to the plate's aspect ratio. 

\item There is a substantial literature on cantilevered beams which bend in two mutually orthogonal directions and also twist about the beam axis.  This corresponds to assuming that there is inextensibility only for deformations along the beam axis and that $v$ is only a function of $x$ and $t$ and $w=h(x,t)+y\alpha (x,t)$, where $x$ is the spatial coordinate along the beam axis (span), $y$ is the spatial coordinate along the chord axis and alpha is the twist about the $x$ axis. For example, see \cite{dowellhodges}. Thus, this is a special case of a nonlinear plate theory (not yet developed here or elsewhere), if only rigid body translation and rotation are considered in the $y$ direction.

\item  In all of our order assumptions herein, we began by declaring $\partial_{x_i} w \sim \eta$. From the physical point of view,  one could make separate assumptions about $w_x$ and $w_y$. An interesting consideration would permit different orders for the slopes of $w$ in the $x$ and $y$ directions, and it would be relevant to further allow different orders in $u$ and $v$.

\end{itemize}

\section{Acknowledgements} The authors acknowledge the generous support of the National Science Foundation in this work through the grant entitled: ``Collaborative Research: Experiment, Theory, and Simulation of Aeroelastic Limit Cycle Oscillations for Energy Harvesting Applications". The first and third authors were partially supported by NSF-DMS 1907620, while the second and fourth authors were partially supported by NSF-DMS 1907500. The second author received substantial support through a Department of Defense SMART Fellowship, and is currently on the research staff of the Air Force Research Laboratory.

\end{document}